\documentclass{amsart}
\usepackage{amsfonts}
\usepackage{amssymb}
\usepackage{amsmath}

\newtheorem{theorem}{Theorem}
\theoremstyle{plain}

\newtheorem{corollary}[theorem]{Corollary}

\newtheorem{remark}{Remark}

\begin{document}
\title[] { Note on the generalization of the higher order $q$-Genocchi numbers and $q$-Euler numbers}
\author{Taekyun Kim  }
\address{Taekyun Kim. Division of General Education-Mathematics \\
Kwangwoon University, Seoul 139-701, Republic of Korea  \\}
\email{tkkim@kw.ac.kr}
\author{Young-Hee Kim}
\address{Young-Hee Kim. Division of General Education-Mathematics \\
Kwangwoon University, Seoul 139-701, Republic of Korea  \\}
\email{yhkim@kw.ac.kr}
\author{Kyung-Won Hwang}
\address{Kyung-Won Hwang. Department of General education\\
Kookmin university, Seoul 136-702, Republic of Korea \\}
\email{khwang7@kookmin.ac.kr}

\maketitle

{\footnotesize {\bf Abstract} \hspace{1mm}
{Cangul-Ozden-Simsek\cite{1} constructed the $q$-Genocchi numbers of
high order using a fermionic $p$-adic integral on $\mathbb{Z}_p$,
and gave Witt's formula and the interpolation functions of these
numbers. In this paper, we present the generalization of the higher
order $q$-Euler numbers and $q$-Genocchi numbers of
Cangul-Ozden-Simsek. We define $q$-extensions of $w$-Euler numbers
and polynomials, and $w$-Genocchi numbers and polynomials of high
order using the multivariate fermionic $p$-adic integral on
$\mathbb{Z}_p$. We have the interpolation functions of these numbers
and polynomials. We obtain the distribution relations for
$q$-extensions of $w$-Euler and $w$-Genocchi polynomials. We also
have the interesting relation for $q$-extensions of these
polynomials. We define $(h,q)$-extensions of $w$-Euler and
$w$-Genocchi polynomials of high order.  We have the interpolation
functions for $(h,q)$-extensions of these polynomials. Moreover, we
obtain some meaningful results of $(h,q)$-extensions of $w$-Euler
and $w$-Genocchi polynomials. }}



\section{Introduction, Definitions and Notations }

Many authors have been studied on the multiple Genocchi and Euler
numbers, and multiple zeta functions (cf. [1-2], [4-6], [9-10],
[14], [17], [19], [22], [24]). In [10], Kim, the first author of
this paper, presented a systematic study of some families of
multiple $q$-Euler numbers and polynomials. By using the
$q$-Volkenborn integration on $\mathbb{Z}_p$, Kim constructed the
$p$-adic $q$-Euler numbers and polynomials of higher order, and gave
the generating function of these numbers and the Euler
$q$-$\zeta$-function. In [14], Kim studied some families of multiple
$q$-Genocchi and $q$-Euler numbers by using the multivariate
$p$-adic $q$-Volkenborn integral on $\mathbb{Z}_p$, and gave
interesting identities related to these numbers.

Recently, Cangul-Ozden-Simsek[1] constructed the $q$-Genocchi
numbers of high order by using a fermionic $p$-adic integral on
$\mathbb{Z}_p$, and gave Witt's formula and the interpolation
functions of these numbers. In [17], Kim gave another constructions
of the $q$-Euler and $q$-Genocchi numbers, which were different from
those of Cangul-Ozden-Simsek. Kim obtained the interesting
relationship between the $q$-$w$-Euler numbers and $q$-$w$-Genocchi
numbers, and gave the interpolation functions of these numbers. In
this paper, we will present the generalization of the higher order
$q$-Euler numbers and $q$-Genocchi numbers of Cangul-Ozden-Simsek
approaching as Kim did in [17].

Throughout this paper, let $p$ be a fixed odd number and the symbols
$\mathbb{Z}_p, \mathbb{Q}_p, \mathbb{C}$ and $\mathbb{C}_p$ denote
the ring of $p$-adic rational integers, the field of $p$-adic
rational numbers, the complex number field and the completion of
algebraic closure of $\mathbb{Q}_p$, respectively. Let $\mathbb{N}$
be the set of natural numbers and $\mathbb{Z}_+ =\mathbb{N}\cup \{0
\}.$ Let $v_p$ be the normalized exponential valuation of
$\mathbb{C}_p$ with $|p|_p =p^{-v_p (p)}= \frac{1}{p}$.

The symbol $q$ can be treated as a complex number, $q \in
\mathbb{C}$, or as a $p$-adic number, $q \in \mathbb{C}_p$. If $q
\in \mathbb{C}$, then we always assume that $|q|< 1.$ If $q \in
\mathbb{C}_p$, then we usually assume that $|1-q|_p< 1$.

Now we will recall some $q$-notations. The $q$-basic natural numbers
are defined by $[n]_q = \frac{1-q^n}{1-q}= 1+q+q^2+ \cdots +
q^{n-1}~ ( n \in \mathbb{N})$, $[n]_{-q} = \frac{1-(-q)^n }{1+q}$
and the $q$-factorial by $[n]_{q}!=[n]_q[n-1]_q \cdots [2]_q[1]_q$.
In this paper, we use the notation $[x]_q = \frac{1-q^x}{1-q}$ and
$[x]_{-q} = \frac{1-(-q)^x}{ 1+q}.$ Hence $ \underset{q \rightarrow
1} {\lim} [x]_q= x$ for any $x$ with $|x|_p \leq 1$ in the present
$p$-adic case (see [1-25]).

The $q$-shift factorial is given by
$$(a:q)_0 =1, \quad (a:q)_k=(1-a)(1-aq) \cdots (1-aq^{k-1}).$$
We note that $\underset{q \rightarrow 1} {\lim}(a:q)_k= (1-a)^k$. It
is known that
$$(a:q)_\infty=(1-a)(1-aq)(1-aq^2)\cdots = \underset {i=1}{\overset{\infty}{\prod}}(1-aq^{i-1}), \quad (\text{see} \,\, [8]).$$
From the definition of the $q$-shift factorial, we note that
$$(a:q)_k=\frac{(a:q)_\infty}{(aq^k:q)_\infty}.$$
Since $\binom{-\alpha}{l}=(-1)^l \binom{\alpha+ l-1}{l}$, it follows
that
\begin{eqnarray*}
\frac{1}{(1-z)^\alpha} &=& (1-z)^{-\alpha}= \underset
{l=0}{\overset{\infty}{\sum}}\binom{-\alpha}{l}(-z)^l= \underset
{l=0}{\overset{\infty}{\sum}}\binom{\alpha+ l-1}{l}z^l.
\end{eqnarray*}

The $q$-binomial theorem is given by
\begin{eqnarray*}\underset{n=0}{\overset{\infty}{\sum}}\frac{(a:q)_n}{(q:q)_n}z^n = \frac{(az:q)_\infty}{(z:q)_\infty},
\end{eqnarray*}
where $z, q \in \mathbb{C}$  with $|z|< 1, \, |q|<1.$ For the
special case, when $a=q^\alpha(\alpha \in \mathbb{C})$, we can write
as follows:
\begin{eqnarray*}
 \frac{1}{(z:q)_\alpha} &=& \frac{(zq^\alpha:q)_\infty}{(z:q)_\infty}
 =\underset{n=0}{\overset{\infty}{\sum}}\frac{(q^\alpha:q)_n}{(q:q)_n}z^n
\\ &=& \underset{n=0}{\overset{\infty}{\sum}}\frac{(1-q^\alpha)(1-q^{\alpha+1})\cdots (1-q^{\alpha+n-1})}
 {(1-q)(1-q^2)\cdots(1-q^n)} z^n\\
 &=&  \underset{n=0}{\overset{\infty}{\sum}}\frac{[\alpha]_q[\alpha+1]_q \cdots [\alpha+n-1]_q}{[1]_q[2]_q \cdots
 [n]_q}z^n \\
 &=& \underset{n=0}{\overset{\infty}{\sum}}\frac{[\alpha]_q[\alpha+1]_q \cdots
 [\alpha+n-1]_q}{[n]_{q}!}z^n .
\end{eqnarray*}
The $q$-binomial coefficients are defined by
\begin{eqnarray*}
\binom{n}{k}_q = \frac{[n]_{q}!}{[k]_{q}![n-k]_{q}!}
=\frac{[n]_q[n-1]_q \cdots [n-k+1]_q}{[k]_q !}, \quad (\text{see}
\,\, [14], [16]).
\end{eqnarray*}
Hence it follows that
$$\frac{1}{(z:q)_\alpha}=\underset{n=0}{\overset{\infty}{\sum}} \binom{n+\alpha-1}{n}_q
z^n,$$ which converges to
$\frac{1}{(1-z)^\alpha}=\underset{n=0}{\overset{\infty}{\sum}}
\binom{n+\alpha-1}{n}z^n$  as $q\rightarrow1.$

We say that $f$ is a uniformly differentiable function at a point $a
\in \mathbb{Z}_p,$ and write $f \in UD (\mathbb{Z}_p),$ the set of
uniformly differentiable function, if the difference quotients $F_g
(x,y)=\frac{f(x)-f(y)}{x-y}$ have a limit $l=f'(a)$ as $(x,
y)\rightarrow (a,a)$. For $f \in UD(\mathbb{Z}_p )$, the
$q$-deformed bosonic $p$-adic integral is defined as
\begin{eqnarray*}
I_q(f)= \int_{\mathbb{Z}_p} f(x)d\mu_q (x)=\underset{N
\rightarrow\infty} {\lim} \underset {x=0}{\overset{p^{N}-1}{\sum}}
f(x)\frac{q^x}{[p^{N}]_q},  \end{eqnarray*} and the $q$-deformed
fermonic $p$-adic integral is defined by
\begin{eqnarray*}
I_{-q} (f)= \int_{\mathbb{Z}_p} f(x)d\mu_{-q} (x)=\underset{N
\rightarrow\infty} {\lim} \underset {x=0}{\overset{p^{N}-1}{\sum}}
f(x)\frac{(-q)^x}{[p^{N}]_{-q}}.
\end{eqnarray*}
The fermionic $p$-adic integral on $\mathbb{Z}_p$ is defined as
$$I_{-1} (f)= \underset{q \rightarrow 1} {\lim}  I_{-q} (f) =\int_{\mathbb{Z}_p} f(x)d\mu_{-1} (x).$$
It follows that $ I_{-1}(f_1)= - I_{-1}(f)+ 2f(0),$ where $f_1 (x)=
f(x+1)$. For details, see [4-17].

The classical Euler polynomials $E_n(x)$ are defined as
\begin{eqnarray*}
\frac{2}{e^t +1}e^{xt}= \underset{x=0}{\overset{\infty}{\sum}} E_n
(x)\frac{t^n}{n!},
\end{eqnarray*}
and the Euler numbers $E_n$ are defined as $E_n= E_n(0)$, (see
[1-25]). The Genocchi numbers are defined as
\begin{equation*}
\frac{2t}{e^{t}+1}=\sum_{n=0}^{\infty}G_{n}\frac{t^{n}}{n!} \quad
{\rm for} \quad |t| < \pi,
\end{equation*}
and the Genocchi polynomials $G_{n}(x)$ are defined as
\begin{equation*}
{\frac{2t}{e^{t}+1}e^{xt}=\sum_{n=0}^{\infty} G_{n}(x)
\frac{t^{n}}{n !}}, \quad  (\text{see  [12], [14], [21]}).
\end{equation*}

It is known that the $w$-Euler polynomials $E_{n,w} (x)$ are defined
as
\begin{eqnarray*}
\frac{2 }{w e^t+1}e^{xt}= \underset{x=0}{\overset{\infty}{\sum}}
E_{n, w} (x)\frac{t^n}{n!},
\end{eqnarray*}
and $E_{n, w} = E_{n, w} (0)$ are called the $w$-Euler numbers. The
$w$-Genocchi polynomials $G_{n, w}(x)$ are defined as
\begin{eqnarray*}
\frac{2t }{w e^t+1}e^{xt}= \underset{x=0}{\overset{\infty}{\sum}}
G_{n, w} (x)\frac{t^n}{n!},
\end{eqnarray*}
and $G_{n, w} = G_{n, w} (0)$ are called the $w$-Genocchi numbers,
(see [1]).

The $w$-Euler polynomials $E_{n,w}^{(r)} (x)$ of order $r$ are
defined as
\begin{eqnarray*}
( \frac{2 }{w e^t+1} )^r e^{xt}=
\underset{x=0}{\overset{\infty}{\sum}} E_{n,w}^{(r)} (x)
\frac{t^n}{n!}, \quad (\text{see} \,\,  [1]),
\end{eqnarray*}
and $E_{n, w}^{(r)} = E_{n, w}^{(r)} (0)$ are called the $w$-Euler
numbers of order $r$. The $w$-Genocchi polynomials $G_{n,
w}^{(r)}(x)$ of order $r$ are defined as
\begin{eqnarray*}
\frac{2t }{w e^t+1}e^{xt}= \underset{x=0}{\overset{\infty}{\sum}}
G_{n, w}^{(r)} (x)\frac{t^n}{n!}, \quad (\text{see} \,\,  [1 ]),
\end{eqnarray*}
and $G_{n, w}^{(r)} = G_{n, w}^{(r)} (0)$ are called the $w$-Euler
numbers of order $r$. When $r=1$ and $w=1$, $E_{n,w}^{(r)} (x)$ and
$E_{n, w}^{(r)}$ are the ordinary Euler polynomials and numbers, and
$G_{n, w}^{(r)}(x)$ and $G_{n, w}^{(r)}$ are the ordinary Genocchi
polynomials and numbers, respectively.

In Section 2, we define $q$-extensions of $w$-Euler numbers and
polynomials of order $r$ and $w$-Genocchi numbers and polynomials of
order $r$, respectively, using the multivariate fermionic $p$-adic
integral on $\mathbb{Z}_p$. We obtain the interpolation functions of
these numbers and polynomials. We have the distribution relations
for $q$-extensions of $w$-Euler polynomials and those of
$w$-Genocchi polynomials. We obtain the interesting relation for
$q$-extensions of these polynomials. We also define
$(h,q)$-extensions of $w$-Euler and $w$-Genocchi polynomials of
order $r$. We have the interpolation functions for
$(h,q)$-extensions of these polynomials. Moreover, we obtain some
meaningful results of $(h,q)$-extensions of $w$-Euler and
$w$-Genocchi polynomials when $h=r-1$.

\vskip 10pt

\section{On the extension of the higher order $q$-Genocchi numbers and $q$-Euler numbers of Cangul-Ozden-Simsek }
\vskip 10pt

In this section, we assume that $w \in \mathbb{C}_p$ with $|1-w|_p <
1$ and $q \in \mathbb{C}_p$ with $|1-q|_p < 1$. Recently,
Cangul-Ozden-Simsek[1] constructed $w$-Genocchi numbers of order
$r$, $G_{n,w}^{(r)}$, as follows.

\begin{eqnarray}
& &t^r \int_{\mathbb{Z}_p^r} w^{x_1 + \cdots+ x_r}e^{t(x_1 + \cdots
+ x_r)}d\mu _{-1} (x_1)\cdots d\mu _{-1} (x_r)  \label{1} \\
& & \qquad = 2^r (\frac{t}{we^t+1})^r
=\underset{n=0}{\overset{\infty}{\sum}}G_{n,w}^{(r)} \frac{t^n}{n!},
\notag \end{eqnarray} where $\int_{\mathbb{Z}_p^r}=
\int_{\mathbb{Z}_{p}} \cdots \int_{\mathbb{Z}_{p}} \, (r-$times) and
$r \in \mathbb{N}$. They also consider the $q$-extension of
$G_{n,w}^{(r)}$ as follows.
\begin{eqnarray}& &t^r\int_{\mathbb{Z}_p^r}
q^{\underset{i=1}{\overset{r}{\sum}}(h-i+1)x_i}e^{t(
\underset{i=1}{\overset{r}{\sum}}x_i)}
d\mu _{-1} (x_1) \cdots d\mu _{-1} (x_r) \label{2}\\
& & \qquad = \frac{2^r t^r}{(q^h e^t+1) \cdots (q^{h-r+1}e^t
+1)}\notag =\underset{n=0}{\overset{\infty}{\sum}}G_{n,q}^{(h,r)}
\frac{t^w}{n!}.\notag
\end{eqnarray}
From $(2)$, they obtained the following interesting formula:
\begin{eqnarray}
G_{n+r,q}^{(r-1,r)}= 2^r r!\binom{n+r}{r} \sum_{v=0}^{\infty}
\binom{r+v-1}{v}_q (-1)^v v^{n}. \label{3} \end{eqnarray}

It seems to be interested to find the numbers corresponding to
$$2^r r!\binom{n+r}{r} \sum_{v=0}^{\infty} \binom{r+v-1}{v}_q (-1)^v
[v]_q^{n}.$$ In the viewpoint of the $q$-extension of $(1)$ using
the multivariate $p-$adic integral on $\mathbb{Z}_{p}$, we define
the $q$-analogue of $w$-Euler numbers of order $r$,
$E_{n,w,q}^{(r)}$, as follows.
\begin{eqnarray}
& &E_{n,w,q}^{(r)}=\int_{\mathbb{Z}_p^r} w^{x_1  + \cdots+ x_r}{[x_1
+ \cdots + x_r]_q ^n}d\mu _{-1} (x_1) \cdots d\mu _{-1} (x_r).
\label{4}\end{eqnarray} From $(\ref{4})$, we note that
\begin{eqnarray*}
E_{n,w,q}^{(r)}&=& \frac{2^r}{(1-q)^n} \sum_{l=0}^{n} \binom{n}{l}
(-1)^l (\frac{1}{1+qw})^r \\
& =&\frac{2^r}{(1-q)^n} \sum_{l=0}^{n} \binom{n}{l} (-1)^l
\sum_{m=0}^{\infty} \binom{m+r-1}{m} (-1)^m q^{lm}w^m \\
& =& 2^r\sum_{m=0}^{\infty} \binom{m+r-1}{m} (-1)^m w^m [m]_q ^n .
\end{eqnarray*}
Therefore, we obtain the following theorem.
\begin{theorem} Let $r \in \mathbb{N}$ and $n \in \mathbb{Z}_+$. Then we have
\begin{eqnarray}E_{n,w,q}^{(r)}=2^r\sum_{m=0}^{\infty} \binom{m+r-1}{m} (-1)^m w^m
[m]_q ^n . \label{5}
\end{eqnarray}
\end{theorem}
\vskip 10pt

Let $F^{(r)}{(t,w|q)}= \underset
{n=0}{\overset{\infty}{\sum}}E_{n,w,q}^{(r)}\frac{t^n}{n!}$. By
$(\ref{4})$ and  $(\ref{5})$, we see that
\begin{eqnarray*}
F^{(r)}{(t,w|q)}&=& \int_{\mathbb{Z}_p^r} w^{x_1+ \cdots+ x_r}e^{
t[x_1 + \cdots + x_r]_q}d\mu _{-1} (x_1) \cdots d\mu _{-1}
(x_r)\\&=&  2^r\sum_{m=0}^{\infty} \binom{m+r-1}{m} (-1)^m w^m
e^{t[m]_q }.
\end{eqnarray*}
Thus we obtain the following corollary.
\begin{corollary}
Let $F^{(r)}{(t,w|q)} = \underset
{n=0}{\overset{\infty}{\sum}}E_{n,w,q}^{(r)}\frac{t^n}{n!}.$ Then we
have
\begin{eqnarray*}
F^{(r)} {(t,w|q)}=2^r\sum_{m=0}^{\infty} \binom{m+r-1}{m} (-1)^m w^m
e^{ t[m]_q}.
\end{eqnarray*}
\end{corollary}
\vskip 10pt

Let us define the $q$-extension of $w$-Euler polynomials of order
$r$ as follows.
\begin{eqnarray}E_{n,w,q}^{(r)}(x)&=& \int_{\mathbb{Z}_p^r} w^{x_1  + \cdots+
x_r}{[x+x_1 + \cdots + x_r]_q ^n}d\mu _{-1} (x_1) \cdots d\mu _{-1}
(x_r).\label{6}
\end{eqnarray}
By $(6)$, we have that
\begin{eqnarray*}
E_{n,w,q}^{(r)}(x) &=& \frac{2^r}{(1-q)^n} \sum_{l=0}^{n}
\binom{n}{l}
(-1)^l q^{lx}(\frac{1}{1+q^lw})^r \\
&=& 2^r\sum_{m=0}^{\infty} \binom{m+r-1}{m} (-1)^m w^m [m+x]_q ^n
\end{eqnarray*}
Therefore, we obtain the following theorem.
\begin{theorem}
Let $r \in \mathbb{N}$ and $n \in \mathbb{Z}_+$. Then we have
\begin{eqnarray}E_{n,w,q}^{(r)}(x)=2^r\underset {m=0}{\overset{\infty}{\sum}} \binom{m+r-1}{m}
(-1)^m w^m [m+x]_q ^n. \label{7}
\end{eqnarray}
\end{theorem}

Let $F^{(r)}{(t,w,x|q)}= \underset
{n=0}{\overset{\infty}{\sum}}E_{n,w,q}^{(r)}(x)\frac{t^n}{n!}.$ By
(\ref{6}) and (\ref{7}), we have
\begin{eqnarray*} F^{(r)}{(t,w,x|q)}
&=& \int_{\mathbb{Z}_p^r} w^{x_1+ \cdots+ x_r}e^{ t[x+x_1 + \cdots +
x_r]_q}d\mu _{-1} (x_1) \cdots d\mu _{-1} (x_r)
\\&=&2^r\underset {m=0}{\overset{\infty}{\sum}} \binom{m+r-1}{m}
(-1)^m w^me^{ t[m+x]_q}. \notag
\end{eqnarray*}
Therefore we have the following corollary.

\begin{corollary}
Let $F^{(r)}{(t,w,x|q)}= \underset
{n=0}{\overset{\infty}{\sum}}E_{n,w,q}^{(r)}(x)\frac{t^n}{n!}.$ Then
we have
\begin{eqnarray} F^{(r)}{(t,w,x|q)}
=2^r\underset {m=0}{\overset{\infty}{\sum}} \binom{m+r-1}{m} (-1)^m
w^me^{ t[m+x]_q}.\label{8}
\end{eqnarray}
\end{corollary}
\vskip 10pt

Now we define the $q$-extension of $w$-Genocchi polynomials of order
$r$, $G_{n,w,q}^{(r)}(x)$, as follows.
\begin{eqnarray}2^rt^r\underset {m=0}{\overset{\infty}{\sum}} \binom{m+r-1}{m} (-1)^m
w^m e^{ t[m+x]_q}=\underset
{n=0}{\overset{\infty}{\sum}}G_{n,w,q}^{(r)}(x)\frac{t^n}{n!}.
\end{eqnarray}
Then we have
\begin{eqnarray}& &
 \quad \underset{n=0}{\overset{\infty}{\sum}}G_{n,w,q}^{(r)}(x)\frac{t^n}{n!} \label{9} \\
&=& t^r\int_{\mathbb{Z}_p^r} w^{x_1+ \cdots+ x_r}e^{t[x+x_1 + \cdots
+ x_r]_q }d\mu _{-1} (x_1) \cdots d\mu _{-1} (x_r) \notag\\
&=&\sum_{n=0}^{\infty}\int_{\mathbb{Z}_{p}^r} w^{x_1+ \cdots+
x_r}{[x+ x_1 + \cdots + x_r]_q ^n}d\mu _{-1} (x_1) \cdots d\mu _{-1}
(x_r) r!\binom{n+r}{r}\frac{t^{n+r}}{(n+r)!}.\notag
\end{eqnarray}
By comparing the coefficients on the both sides of $(\ref{9})$, we
see that $$G_{0,w,q}^{(r)}(x)=G_{1,w,q}^{(r)}(x)= \cdots
=G_{r-1,w,q}^{(r)}(x)=0, $$ and
\begin{eqnarray} & & G_{n+r,w,q}^{(r)}(x) \label{11}\\
&=& r!\binom{n+r}{r} \int_{\mathbb{Z}_p^r} w^{x_1 +x_2 + \cdots+
x_r}{[x+x_1 + \cdots +
x_r]_q ^n}d\mu _{-1} (x_1) \cdots d\mu _{-1} (x_r) \notag \\
&=& r!\binom{n+r}{r}E_{n,w,q}^{(r)}(x).\notag
\end{eqnarray} In the
special case of $x=0$, $G_{n,w,q}^{(r)}(0)=G_{n,w,q}^{(r)}$ are
called the $q$-extension of $w$-Genocchi numbers of order $r$. By
$(\ref{11})$, we have the following theorem.

\begin{theorem}
Let $r \in \mathbb{N}$ and $n \in \mathbb{Z}_+$. Then we have
\begin{eqnarray*}\frac{G_{n+r,w,q}^{(r)}(x)}{r!\binom{n+r}{r}}&=& \int_{\mathbb{Z}_p^r} w^{x_1+ \cdots+
x_r}{[x+x_1 + \cdots + x_r]_q ^n}d\mu _{-1} (x_1) \cdots d\mu _{-1}
(x_r)\\
&=&E_{n,w,q}^{(r)}(x),
\end{eqnarray*}
and $\, G_{0,w,q}^{(r)}(x)=G_{1,w,q}^{(r)}(x)= \cdots
=G_{r-1,w,q}^{(r)}(x)=0.$
\end{theorem}
\vskip 10pt

Now we consider the distribution relation for the $q$-extension of
$w$-Euler polynomials of order $r$. For $d \in \mathbb{N}$ with
$d\equiv 1$ (mod $2$), by $(\ref{8})$, we see that

\begin{eqnarray}
&& \qquad F^{(r)}{(t,w,x|q)} \label{12}\\
&= &2^r\underset {m=0}{\overset{\infty}{\sum}}
\binom{m+r-1}{m} (-1)^m w^m e^{ t[m+x]_q} \notag \\
&=& \sum_{a_1, \cdots a_r=0}^{d-1}(\underset
{i=1}{\overset{r}{\prod}}w^{a_i}) (-1)^{a_1 + \cdots + a_r}2^r
\underset {m=0}{\overset{\infty}{\sum}} \binom{m+r-1}{m} (-1)^m
w^{md} e^{ t [d]_q [m+\frac{a_1 + \cdots + a_r +x}{d}]_{q^d} }  \notag\\
&=&\sum_{a_1, \cdots a_r=0}^{d-1}(\underset
{i=1}{\overset{r}{\prod}}w^{a_i}) (-1)^{a_1 + \cdots + a_r}F
^{(r)}{([d]_q t,\, w^d, \, \frac{a_1 + \cdots + a_r +x}{d}\, | \,
q^d)}.\notag
\end{eqnarray}
By $(\ref{12})$, we obtain the following distribution relations for
$E_{n,w,q}^{(r)}(x)$ and $G_{n+r,w,q}^{(r)}(x)$, respectively.
\begin{theorem}
Let $r \in \mathbb{N}$, $n \in \mathbb{Z}_+$ and $d \in \mathbb{N}$
with $d \equiv 1$ (mod $2$). Then we have
\begin{eqnarray*}E_{n,w,q}^{(r)}(x) = [d]_q ^n \sum_{a_1, \cdots a_r=0}^{d-1}(\underset
{i=1}{\overset{r}{\prod}}w^{a_i}) (-1)^{a_1 + \cdots + a_r}
E_{n,w^d,q^d}^{(r)}(\frac{a_1 + \cdots + a_r+x}{d}) \notag .
\end{eqnarray*}
Furthermore,
\begin{eqnarray*} G_{n+r,w,q}^{(r)}(x) = [d]_q ^n \sum_{a_1, \cdots a_r=0}^{d-1}(\underset
{i=1}{\overset{r}{\prod}}w^{a_i}) (-1)^{a_1 + \cdots + a_r}
G_{n+r,w^d,q^d}^{(r)}(\frac{a_1 + \cdots + a_r+x}{d}) \notag.
\end{eqnarray*}
\end{theorem}


\vskip 10pt

For the extension of $(\ref{2})$, we consider the $(h,q)$-extension
of $w$-Euler polynomials of order $r$. For $h \in \mathbb{Z}$, $r
\in \mathbb{N}$ and $n \in \mathbb{Z}_+$, let us define the
$(h,q)$-extension of $w$-Euler polynomial of order $r$ as follows.

\begin{eqnarray} & &E_{n,w,q}^{(h,r)}(x) \label{13}\\
& &= \int_{\mathbb{Z}_p^r} w^{x_1+ \cdots+ x_r}{[x+x_1 + \cdots +
x_r]_q ^n} q^{\underset{i=1}{\overset{r}{\sum}}(h-i+1)x_i}d\mu _{-1}
(x_1) \cdots d\mu _{-1} (x_r).\notag
\end{eqnarray}
From $(\ref{13})$, we note that
\begin{eqnarray}
E_{n,w,q}^{(h,r)}(x) &=& \frac{2^r}{(1-q)^n}
\underset{l=0}{\overset{n}{\sum}}\frac{\binom{n}{l} (-1)^l
q^{lx}}{(1+q^{l+h}w)(1+q^{l+h-1}w)\cdots (1+q^{l+h-r+1}w)} \notag \\
&=& \frac{2^r}{(1-q)^n}
\underset{l=0}{\overset{n}{\sum}}\frac{\binom{n}{l} (-1)^l
q^{lx}}{(-q^{l+h}w: q^{-1})_r} \label{14}\\
&=& \frac{2^r}{(1-q)^n}
\underset{l=0}{\overset{n}{\sum}}\binom{n}{l} (-1)^l
q^{lx}\underset{m=0}{\overset{\infty}{\sum}}{\binom{m+r-1}{m}}_{q^{-1}}(-1)^m
q^{(l+h)m}w^m \notag \\
&=&2^r\underset{m=0}{\overset{\infty}{\sum}}{\binom{m+r-1}{m}}_{q^{-1}}(-1)^m
q ^{hm}w^m [m+x]_q ^n. \notag
\end{eqnarray}
Therefore, we obtain the following theorem.
\begin{theorem}
Let $h \in \mathbb{Z}$, $r \in \mathbb{N}$ and $n \in \mathbb{Z}_+$.
Then we have
\begin{eqnarray}
E_{n,w,q}^{(h,r)}(x) &=& \frac{2^r}{(1-q)^n}
\underset{l=0}{\overset{n}{\sum}}\frac{\binom{n}{l} (-1)^l
q^{lx}}{(-q^{l+h}w: q^{-1})_r}  \\
&=&2^r\underset{m=0}{\overset{\infty}{\sum}}{\binom{m+r-1}{m}}_{q^{-1}}(-1)^m
q ^{hm}w^m [m+x]_q ^n.\notag
\end{eqnarray}
\end{theorem}
We also have the following result.

\begin{corollary}
Let $F^{(h,r)}{(t,w,x|q)}= \underset
{n=0}{\overset{\infty}{\sum}}E_{n,w,q}^{(h,r)}(x)\frac{t^n}{n!}.$
Then we have
\begin{eqnarray}F^{(h,r)}{(t,w,x|q)}=2^r\underset{m=0}{\overset{\infty}{\sum}}{\binom{m+r-1}{m}}_{q^{-1}}(-1)^m
q ^{hm}w^m e^{ t[m+x]_q} .
\end{eqnarray}
\end{corollary}

\begin{remark}
In the special case $x=0$, $E_{n,w,q}^{(h,r)}(0)= E_{n,w,q}^{(h,r)}$
are called the $(h,q)$-extension of $w$-Euler numbers of order $r$.
\end{remark}

\vskip 10pt

If we take $h=r-1$ in $(\ref{14})$, then we have
\begin{eqnarray}
E_{n,w,q}^{(r-1,r)}(x) &=& \frac{2^r}{(1-q)^n}
\underset{l=0}{\overset{n}{\sum}}\frac{\binom{n}{l} (-1)^l
q^{lx}}{(1+q^{l+r-1}w)(1+q^{l+r-2}w)\cdots (1+q^{l}w)} \notag \\
&=& \frac{2^r}{(1-q)^n}
\underset{l=0}{\overset{n}{\sum}}\frac{\binom{n}{l} (-1)^l
q^{lx}}{(-q^{l}w: q)_r} \\
&=& \frac{2^r}{(1-q)^n}
\underset{l=0}{\overset{n}{\sum}}\binom{n}{l} (-1)^l
q^{lx}\underset{m=0}{\overset{\infty}{\sum}}{\binom{m+r-1}{m}}_{q}(-1)^m
q^{lm}w^m \notag \\
&=&2^r\underset{m=0}{\overset{\infty}{\sum}}{\binom{m+r-1}{m}}_{q}(-1)^m
w^m [m+x]_q ^n. \notag
\end{eqnarray}
Then we have the following theorem.
\begin{theorem} Let $r \in \mathbb{N}$ and $n \in \mathbb{Z}_+$.
Then we have
\begin{eqnarray*}
E_{n,w,q}^{(r-1,r)}(x) &=& \frac{2^r}{(1-q)^n}
\underset{l=0}{\overset{n}{\sum}}\frac{\binom{n}{l} (-1)^l
q^{lx}}{(-q^{l}w: q)_r}\\
&=&2^r\underset{m=0}{\overset{\infty}{\sum}}{\binom{m+r-1}{m}}_{q}(-1)^m
w^m [m+x]_q ^n.
\end{eqnarray*}
\end{theorem}

We also the following corollary.
\begin{corollary}
Let $F^{(r-1,r)}{(t,w,x|q)}= \underset
{n=0}{\overset{\infty}{\sum}}E_{n,w,q}^{(r-1,r)}(x)\frac{t^n}{n!}.$
Then we have
\begin{eqnarray}
F^{(r-1,r)}{(t,w,x|q)} =
2^r\underset{m=0}{\overset{\infty}{\sum}}{\binom{m+r-1}{m}}_{q}(-1)^m
w^m e^{t[m+x]_q }. \label{18}
\end{eqnarray}
\end{corollary}

From $(\ref{18})$, we note that

\begin{eqnarray}
F^{(r-1,r)}{(t,w,x|q)}&=&
2^r\underset{m=0}{\overset{\infty}{\sum}}{\binom{m+r-1}{m}}_{q}(-1)^m
w^m e^{ t[m+x]_q} \notag \\
&=& \sum_{a_1, \cdots a_r=0}^{d-1}q^{\underset
{i=0}{\overset{r}{\sum}}(r-i)a_i} (-1)^{a_1 + \cdots +
a_r}w^{a_1 + \cdots + a_r} \label{19} \\
&& \quad \times 2^r \underset {m=0}{\overset{\infty}{\sum}}
{\binom{m+r-1}{m}}_{q^d}(-1)^m
w^{md}e^{ t [d]_q [m+\frac{a_1 + \cdots + a_r +x}{d}]_{q^d} } \notag \\
&=& \sum_{a_1, \cdots a_r=0}^{d-1}q^{\underset
{i=0}{\overset{r}{\sum}}(r-i)a_i} (-1)^{a_1 + \cdots + a_r}w^{a_1 +
\cdots + a_r} \notag \\
&&  \quad \times F^{(r-1,r)}{([d]_q t,w^d,\frac{a_1 + \cdots + a_r
+x}{d}|q^d)},\notag
\end{eqnarray}
where $d \in \mathbb{N}$ with $d\equiv 1$ (mod  $2$). By
$(\ref{19})$, we obtain the following the distribution relation for
$E_{n,w,q}^{(r-1,r)}(x)$.

\begin{theorem} For  $r \in
\mathbb{N}$, $n \in \mathbb{Z}_+$ and  $d \in \mathbb{N}$ with
$d\equiv 1 $ (mod $2$). Then we have
\begin{eqnarray*}
&&E_{n,w,q}^{(r-1,r)}(x) \\
&&= [d]_q ^n \sum_{a_1, \cdots a_r=0}^{d-1}q^{\underset
{i=0}{\overset{\infty}{\sum}}(r-i)a_i} (-1)^{a_1 + \cdots +
a_r}w^{a_1 + \cdots + a_r}E_{n,w^d,q^d}^{(r-1,r)}(\frac{a_1 + \cdots
+ a_r +x}{d}).\notag
\end{eqnarray*}

\end{theorem}
\vskip 10pt

Now we define the $(h,q)$-extension of $w$-Genocchi polynomials
$G_{n,w,q}^{(h,r)}(x)$ of order $r$ as follows.
\begin{eqnarray}
2^rt^r\underset{m=0}{\overset{\infty}{\sum}}{\binom{m+r-1}{m}}_{q^{-1}}(-1)^m
q^{hm}w^m e^{ t[m+x]_q} = \underset
{n=0}{\overset{\infty}{\sum}}G_{n,w,q}^{(h,r)}(x)\frac{t^n}{n!}.
\end{eqnarray}
Then we have
\begin{eqnarray}
& &\underset
{n=0}{\overset{\infty}{\sum}}G_{n,w,q}^{(h,r)}(x)\frac{t^n}{n!}\notag\\
& & \quad=t^r \int_{\mathbb{Z}_p^r} q^{\underset
{i=0}{\overset{\infty}{\sum}}(h-i+1)x_i}w^{x_1+ \cdots+ x_r}e^{
t[x+x_1 + \cdots + x_r]_q}d\mu _{-1} (x_1) \cdots d\mu _{-1}
(x_r) \label{21}\\
& & \quad = \underset
{n=0}{\overset{\infty}{\sum}}\int_{\mathbb{Z}_p^r} q^{\underset
{i=0}{\overset{\infty}{\sum}}(h-i+1)x_i}w^{x_1+ \cdots+ x_r}{[x+x_1
+ \cdots + x_r]_q ^n }d\mu _{-1} (x_1) \cdots
d\mu _{-1} (x_r) \notag\\
&&\qquad \times r!\binom{n+r}{r}\frac{t^{n+r}}{(n+r)!}.\notag
\end{eqnarray}
From $(\ref{13})$ and $(\ref{21})$, we derive the following result.

\begin{theorem} Let $r \in \mathbb{N}$ and $n \in \mathbb{Z}_+$.
Then we have
\begin{eqnarray*}
\frac{G_{n+r,w,q}^{(h,r)}(x)}{r!\binom{n+r}{r}}&=&
\int_{\mathbb{Z}_p^r} q^{\underset
{i=0}{\overset{\infty}{\sum}}(h-i+1)x_i}w^{x_1+ \cdots+ x_r}{[x+x_1
+ \cdots + x_r]_q ^n }d\mu _{-1} (x_1) \cdots d\mu _{-1}
(x_r)  \\
&= &E_{n,w,q}^{(h,r)}(x),\notag
\end{eqnarray*} and
$\, G_{0,w,q}^{(h,r)}(x)=G_{1,w,q}^{(h,r)}(x)= \cdots
=G_{r-1,w,q}^{(h,r)}(x)=0$.
\end{theorem}

\vskip 10pt

When $h=r-1$ in Theorem 12, we have
\begin{eqnarray*}
\frac{G_{n+r,w,q}^{(r-1,r)}(x)}{r!\binom{n+r}{r}}&=&
\int_{\mathbb{Z}_p^r} q^{\underset
{i=0}{\overset{\infty}{\sum}}(r-i)x_i}w^{x_1+ \cdots+ x_r}{[x+x_1 +
\cdots + x_r]_q ^n }d\mu _{-1} (x_1) \cdots d\mu _{-1}
(x_r)  \\
&=&
2^r\underset{m=0}{\overset{\infty}{\sum}}{\binom{m+r-1}{m}}_{q}(-1)^m
w^m [m+x]_q ^n \\
&=& \frac{2^r}{(1-q)^n}
\underset{l=0}{\overset{n}{\sum}}\frac{\binom{n}{l} (-1)^l
q^{lx}}{(-q^{l}w: q)_r} \\
&=& E_{n,w,q}^{(r-1,r)}(x).\end{eqnarray*}

\begin{remark}
In the special case $x=0$, $G_{n,w,q}^{(h,r)}(0)= G_{n,w,q}^{(h,r)}$
are called the $(h,q)$-extension of $w$-Genocchi numbers of order
$r$.
\end{remark}


\bigskip

\begin{thebibliography}{99}
\bibitem{1} I. N. Cangul, H. Ozden, Y. Simsek,  \textit{A new approach to $q$-Genocchi numbers
and their interpolations}, Nonlinear Analysis (2008),
doi:10.1016/j.na.2008.11.040.

\bibitem{2} M. Cenkci, Y. Simsek, V. Kurt, \textit{Multiple two-variable $p$-adic
$q$-$L$-function and its behavior at $s = 0$},  Russ. J. Math. Phys.
\textbf{15} (2008) no. 4, 447-459.


\bibitem{3} G. Kim, B. Kim, J. Choi, \textit{The DC
algorithm for computing sums of powers of consecutive integers and
Bernoulli numbers}, Adv. Stud. Contemp. Math. (Kyungshang)
\textbf{17} (2008), no. 2, 137--145.


\bibitem{4} T. Kim, \textit{$q$-Volkenborn integration},
Russ. J. Math. Phys.  \textbf{9} (2002), 288--299.

\bibitem{5}  T. Kim, \textit{On Euler-Barnes multiple zeta
functions}, Russ. J. Math. Phys. \textbf{10} (2003), no. 3,
261--267.

\bibitem{6}  T. Kim, \textit{Analytic continuation of
multiple $q$-zeta functions and their values at negative integers},
Russ. J. Math. Phys. \textbf{11} (2004), no. 12, 71--76.

\bibitem{7}  T. Kim, \textit{Power series and asymptotic
series associated with the $q$-analog of the two-variable $p$-adic
$L$-function}, Russ. J. Math. Phys. \textbf{12} (2005), no. 2,
186--196.

\bibitem{8} T. Kim, \textit{$q$-Generalized Euler numbers and polynomials},
Russ. J. Math. Phys.  \textbf{13} (2006), no. 3, 293--298.

\bibitem{9}  T. Kim, \textit{Multiple $p$-adic $L$-function},
Russ. J. Math. Phys. \textbf{13} (2006), no. 2, 151--157.

\bibitem{10} T. Kim, \textit{$q$-Euler numbers and polynomials associated with
$p$-adic $q$-integrals}, J. Nonlinear Math. Phys. \textbf{14}
(2007), no. 1, 15--27.

\bibitem{11} T. Kim, \textit{A note on $p$-adic $q$-integral on $\mathbb{Z}_p$
associated with $q$-Euler numbers}, Adv. Stud. Contemp. Math.
(Kyungshang) \textbf{15} (2007), 133--138.

\bibitem{12} T. Kim, \textit{On the $q$-extension of Euler and Genocchi numbers},
J. Math. Anal. Appl. \textbf{326} (2007), 1458--1465.

\bibitem{13}  T. Kim, \textit{$q$-Extension of the Euler
formula and trigonometric functions}, Russ. J. Math. Phys.
\textbf{14} (2007), no. 3, 275--278.


\bibitem{14} T. Kim, \textit{On the multiple $q$-Genocchi and Euler numbers},
Russ. J. Math. Phys.  \textbf{15} (2008) no.4, 481-486.

\bibitem{15} T. Kim, \textit{The modified $q$-Euler numbers and
polynomials}, Adv. Stud. Contemp. Math. (Kyungshang) \textbf{16}
(2008), 161--170.

\bibitem{16} T. Kim, \textit{$q$-Bernoulli numbers and polynomials associated
with Gaussian binomial coefficients}, Russ. J. Math. Phys.
\textbf{15} (2008), no.1, 51-57.

\bibitem{17} T. Kim, \textit{New approach to $q$-Euler, Genocchi numbers and their interpolation functions},
 arXiv:0901.0353v1 [math.NT].



\bibitem{18} Y.-H. Kim, W. Kim, L.-C. Jang, \textit{On the $q$-extension of
Apostol-Euler numbers and polynomials}, Abstr. Appl. Anal.
\textbf{2008} (2008), Article ID 296159, 10 pages.


\bibitem{19} H. Ozden, I. N. Cangul, Y. Simsek, \textit{Multivariate interpolation
functions of higher-order $q$-Euler numbers and their applications},
Abstr. Appl. Anal. \textbf{2008} (2008), Art. ID 390857, 16 pages.


\bibitem{20} H. Ozden, Y. Simsek, \textit{A new extension of $q$-Euler numbers and polynomials
related to their interpolation functions}, Appl. Math. Lett.
\textbf{21} (2008), 934--939.


\bibitem{21} K. H. Park, Y.-H. Kim, \textit{On some arithmetical properties of the
Genocchi numbers and polynomials}, Advances in Difference Equations
(2008), http://www.hindawi.com/journals/ ade/aip.195049.html.

\bibitem{22} Y. Simsek, \textit{Complete sum of products of $(h,q)$-extension of the Euler polynomials and numbers},
 arXiv:0707.2849v1 [math.NT].

\bibitem{23} J. V. Leyendekkers, A. G. Shannon, C. K. Wong,  \textit{Integer structure
analysis of the product of adjacent integers and Euler's extension
of Fermat's last theorem}, Adv. Stud. Contemp. Math. (Kyungshang)
\textbf{17} (2008), no. 2, 221--229.

\bibitem{24} H.M. Srivastava, T. Kim, Y. Simsek,
\textit{$q$-Bernoulli numbers and polynomials associated with
multiple $q$-zeta functions and basic $L$-series}, Russ. J. Math.
Phys. \textbf{12} (2005), no. 2, 241--268.

\bibitem{25} Z. Zhang, Y. Zhang, \textit{Summation formulas of $q$-series by
modified Abel's lemma}, Adv. Stud. Contemp. Math. (Kyungshang)
\textbf{17} (2008), no. 2, 119--129.




\end{thebibliography}
\end{document}